\documentclass[10pt,a4paper,onecolumn,lmr]{article}
\usepackage[utf8]{inputenc}
\usepackage[english]{babel}
\usepackage{amsmath}
\usepackage{amsfonts}
\usepackage{amssymb}
\usepackage{setspace}
\usepackage{amsthm}
\usepackage{blindtext}
\usepackage[numbers]{natbib}
\theoremstyle{plain}
\newtheorem{theorem}{Theorem}[section]

\newtheorem{lemma}[theorem]{Lemma}

\theoremstyle{definition}
\newtheorem{definition}[theorem]{Definition} 
\newtheorem{example}[theorem]{Example}
\newtheorem{assumption}[]{Assumption}

\title{Ruin Probabilities in a \\ Markovian Shot-Noise Environment}
\author{Simon Pojer, Stefan Thonhauser \\ Graz University of Technology }
\date{August 2021  }
\begin{document}
\onehalfspacing
\maketitle
\begin{abstract}
We consider a risk model with a counting process whose intensity is a Markovian shot-noise process, to resolve one of the disadvantages of the Cramér-Lundberg model, namely the constant jump intensity of the Poisson process. Due to this structure, we can apply the theory of PDMPs on a multivariate process containing the intensity and the reserve process, which allows us to identify a family of martingales. Eventually, we use change of measure techniques to derive an upper bound for the ruin probability in this model. Exploiting a recurrent structure of the shot-noise process, even the asymptotic behaviour of the ruin probability can be determined.
\end{abstract} 

\section{Introduction}
The theory of doubly stochastic Poisson processes described in \cite{Bremaud.1981}, allows the generalization of the well-known Cramér-Lundberg model to the broad class of Cox models, which are discussed e.g. in \cite{Grandell.1991}.
Members of this family are for example the Markov-modulated risk model, where the intensity is modelled by a continuous-time Markov chain which can be found in \cite[][Chapter VII]{Asmussen.2010} and \cite[][Chapter 8]{Rolski.1999}, the Björk-Grandell model considered in \cite{Schmidli.1997} and diffusion-driven models studied in \cite{Grandell.2011}. \\ 
Especially, arrivals of claims caused by catastrophic events can be realistically modelled using shot-noise intensity. This has been done in \cite{Albrecher.2006}, \cite{Macci.2011} and \cite{Dassios.2015} where the asymptotic behaviour of the ruin probability in general shot-noise environments was studied. In these settings, upper and lower bounds could be derived. The idea of applying the theory of piecewise deterministic Markov processes to a Cox model with Markovian shot-noise intensity was used in \cite{Dassios.2003} and \cite{Dassios.2005} in the context of pricing reinsurance contracts.\\ \\
 Interested in the behaviour of the ruin probability in this model, we follow the PDMP approach to find suitable alternative probability measures. Further, we take advantage of the properties of the process under these measures to obtain an exponential decreasing upper bound. Exploiting a recurrent behaviour of the shot-noise process and applying the extended renewal theory obtained in \cite{Schmidli.1997} we eventually derive the exact asymptotic behaviour of the ruin probability. 

\section{The Markovian Shot-Noise Ruin Model}
We assume for the rest of this paper the existence of a complete probability space $\left(\Omega, \mathcal{F}, \mathbb{P}\right)$ which is big enough to contain all mentioned stochastic processes and random variables. For some stochastic process $Z$ we denote $\mathcal{F}^Z_t$ for the $\mathbb{P}$-complete and right continuous natural filtration.
Now let us define the shot-noise environment properly. For this we consider the following four objects: A Poisson process $N^\lambda$ with constant intensity $\rho>0$ and jump times $\left\lbrace T^\lambda_i\right\rbrace_{i \in \mathbb{N}}$, a sequence $\left\lbrace Y_i\right\rbrace_{i \in \mathbb{N}}$ of positive i.i.d. random variables with distribution function $F_Y$, a non-negative function $w$, and a positive starting value $\lambda_0$. With these components we define the multiplicative shot-noise process by \[ \lambda_t:= \lambda_0w(t) + \sum_{i=1}^{N^\lambda_t}Y_i w(t-T^\lambda_i).\] 
Since we want to exploit the theory of PDMPs it would be preferable if the process $\lambda$ satisfies the Markov property. This is equivalent to the existence of some $\delta>0$ such that $w(t) = e^{-\delta t}$. Due to this, we define the Markovian shot-noise process the following way. 

\begin{definition}
Let $N^\lambda$ be a Poisson process with intensity $\rho >0$ and jump times $\left\lbrace T^\lambda_i\right\rbrace_{i \in \mathbb{N}}$, $\left \lbrace Y_i \right \rbrace _{i \in \mathbb{N}}$ i.i.d. copies of a positive random variable $Y$ with distribution function $F_Y$ and independent of the process $N^\lambda$, $\lambda_0>0$ and $\delta >0$ constant. Then we define the Markovian shot-noise process as
\[ \lambda_t = \lambda_0e^{-\delta t} + \sum_{i=1}^{N^\lambda_t}Y_i e^{-\delta(t-T^\lambda_i)}.\] 
\end{definition}
As shown in \cite{Dassios.2005} the Markovian shot-noise process is a piecewise deterministic Markov process with generator \[ \mathcal{A}^\lambda f(\lambda) = - \delta \lambda \frac{\partial f(\lambda)}{\partial \lambda} + \rho \int_0^\infty \left(f(\lambda + y) - f(\lambda)\right) \, F_Y(\mathrm{d}y). \] Further information about PDMPs can be found in \cite{Davis.1993}  or \cite[Chapter 11]{Rolski.1999}.
To fully specify our model we will now define the surplus process.
\begin{definition}
Let $\lambda$ be a Markovian shot-noise process, $N$ a Cox process with intensity $\lambda$ and a sequence $\left\lbrace U_i\right\rbrace_{i \in \mathbb{N}}$ of i.i.d. copies  of a positive random variable $U$ with continuous distribution $F_U$ which are independent of $N$. For some initial capital $u$ and constant premium rate $c>0$ we define the surplus process by \[ X_t= u + ct - \sum_{i=1}^{N_t} U_i.\]
\end{definition}
Now define $\mathcal{F}_t:= \mathcal{F}^X_t \vee \mathcal{F}^\lambda_t$, hence $\left\lbrace\mathcal{F}_t\right\rbrace_{t \geq 0}$ is the combined filtration of the Markovian shot-noise process and the surplus process. If not mentioned differently, we will from now on consider the filtered probability space $\left(\Omega, \mathcal{F}, \left\lbrace \mathcal{F}_t\right\rbrace_{t \geq 0}, \mathbb{P}_{(u,\lambda_0)}\right),$ where we define the measure $\mathbb{P}_{(u,\lambda_0)}$ as the measure $\mathbb{P}$ under the conditions that the initial capital of the surplus process is $u$ and the starting intensity is $\lambda_0$. We will denote the expectation of a random variable $Z$ under this measure by $\mathbb{E}_{(u,\lambda_0)}\left[Z\right]$ or $\mathbb{E}\left[Z\right]$ if $Z$ is independent of the initial values. \\ \\
The multivariate process $\left(X, \lambda, \cdot\right):= \left((X_t,\lambda_t,t)\right)_{t \geq 0}$ is a càdlàg PDMP without active boundary and generator \begin{align*}
\mathcal{A}f(x,\lambda,t) &= c\frac{\partial f(x,\lambda, t)}{\partial x} - \delta \lambda\frac{\partial f(x,\lambda, t)}{\partial \lambda} + \frac{\partial f(x,\lambda, t)}{\partial t}\\[7pt]
& + \lambda \int_0^\infty \left(f(x-u,\lambda,t)-f(x,\lambda,t)\right) \, F_U(\mathrm{d}u) \\[7pt]
&+ \rho \int_0^\infty \left(f(x,\lambda+y,t)-f(x,\lambda,t)\right)\, F_Y(\mathrm{d}y).
\end{align*}
Its domain consists of all functions $f$ which are absolutely continuous and satisfy the integrability condition \[\mathbb{E}_{(u,\lambda_0)} \left[ \sum_{i=1}^{\tilde N_t} \left\vert f(X_{T_i}, \lambda_{T_i}, T_i) - f(X_{T_i-}, \lambda_{T_i-}, T_i-)\right\vert \right] < \infty\] for all $t\geq 0$, where $\tilde N$ denotes the process counting the random jumps of the PDMP $(X,\lambda,\cdot)$. Similar to the Cramér-Lundberg model we want to state a net profit condition, which is necessary to ensure that ruin does not occur with probability $1$.
\begin{lemma}\label{lemma_expectation_surplus}
The surplus process satisfies \[ \lim_{t \to \infty} \frac{\mathbb{E}_{(u,\lambda_0)}\left[X_t\right]}{t} = c- \frac{\rho}{\delta} \mathbb{E} \left[U\right]\mathbb{E}\left[Y\right].\]
\begin{proof}
The function $\bar f(x,\lambda,t):= x$ is in the domain of the generator. Consequently, \begin{align*}
\mathbb{E}_{(u,\lambda_0)}\left[X_t\right] &= u + \mathbb{E}_{(u,\lambda_0)}\left[\int_0^t \mathcal{A}\bar f(X_s,\lambda_s,s) \, \mathrm{d}s\right] = u +ct - \mathbb{E}_{(u,\lambda_0)} \left[ \int_0^t \lambda_s \mathbb{E}\left[U\right] \, \mathrm{d}s\right].
\end{align*}
The process $\lambda$ is positive so we can use Tonelli's theorem and interchange expectation and integration, which leads to 
\begin{align} \label{equation_expectation_X}
\mathbb{E}_{(u,\lambda_0)}\left[X_t\right] = u +ct- \mathbb{E}\left[U\right]\int_0^t \mathbb{E}_{(u,\lambda_0)} \left[\lambda_s\right] \, \mathrm{d}s.
\end{align}

Now we use the same procedure to obtain an equation for $\mathbb{E}_{(u,\lambda_0)}\left[\lambda_s\right]$. Defining the function $\tilde f(x,\lambda,t):= \lambda$ we get \begin{align*}
\mathbb{E}_{(u,\lambda_0)}\left[\lambda_s\right] = \lambda_0 - \delta \int_0^s \mathbb{E}_{(u,\lambda_0)}\left[\lambda_u\right]\, \mathrm{d}u + \rho s \mathbb{E}\left[Y\right].
\end{align*}
Differentiating both sides with respect to $s$ gives us that $\mathbb{E}_{(u,\lambda_0)}\left[\lambda_s\right]$ is the solution to the differential equation $ g'(s) = -\delta g(s) + \rho \, \mathbb{E}\left[Y\right],$ with initial value $g(0)=\lambda_0$. Solving the ODE gives us \begin{align}\label{equation_expectation_lambda}
 \mathbb{E}_{(u,\lambda_0)}\left[\lambda_s\right] = \lambda_0e^{-\delta s} + \frac{\rho}{\delta} 
\mathbb{E}\left[Y\right] (1-e^{-\delta s}).
\end{align} 
Using the result \ref{equation_expectation_lambda} in equation \ref{equation_expectation_X} leads to \begin{align*}
\mathbb{E}_{(u,\lambda_0)}\left[X_t\right] = u &+ ct- 
\mathbb{E}\left[U\right]\frac{\rho}{\delta}  \mathbb{E}\left[Y\right] t \\ &+ \mathbb{E}\left[U\right] \left(\frac{\lambda_0}{\delta} - \frac{\rho}{\delta^2} \mathbb{E}\left[Y\right] \right)\left(1-e^{-\delta t}\right).
\end{align*} 
Now, let us divide by $t$ and let it tend to infinity to obtain \[  \lim_{t \to \infty} \frac{\mathbb{E}_{(u,\lambda_0)}\left[X_t\right]}{t} = c- \frac{\rho}{\delta} \mathbb{E}\left[U\right]\mathbb{E}\left[Y\right].\]
\end{proof}
\end{lemma}
Motivated by this result we assume the following:
\begin{assumption}
We assume that the net profit condition \[c> \frac{\rho}{\delta} \mathbb{E} \left[U\right]\mathbb{E}\left[Y\right], \] is satisfied.
\end{assumption}

\section{Martingales and Change of Measure}
To obtain the asymptotic behaviour of the ruin probability in this model, we want to exploit the following result derived in \cite{Schmidli.1997}.

\begin{theorem}{ \cite[][Theorem 2]{Schmidli.1997}}\label{renewal_theorem_schmidli}\\
Assume that $z(u)$ is directly Riemann integrable, that $0\leq p(u,x) \leq 1$ is continuous in $u$ and that$\int_0^u \!p(u,y) \, B(\mathrm{d}y)$ is directly Riemann integrable. Denote by $Z(u)$ the solution to \[ Z(u) = \int_0^u \!Z(u-y)(1-p(u,y)) \, B(\mathrm{d}y) + z(u),\] which is bounded on bounded intervals. Then the limit \[ \lim_{u \to \infty} Z(u) \] exists and is finite provided $B(u)$ is not arithmetic. If $B(u)$ is arithmetic with span $\gamma$, then \[ \lim_{n \to \infty} Z(x+n\gamma) \] exists and is finite for all $x$ fixed.
\end{theorem}
Unfortunately, we cannot apply this Theorem directly to the ruin probability because of two problems. The first issue is, that the ruin probability depends on the initial intensity level $\lambda_0$. To bypass this, we have to choose appropriate renewal times such that $\lambda$ has always the same level, which we will do in Section \ref{section_renewal_equation}. The second problem is, that suitable choices of $B$ are defective under the original measure $\mathbb{P}_{(u,\lambda_0)}$. This is a common issue and can be solved through change of measure techniques. \\ \\
To do so we have to find martingales of the form $M_t=h(X_t, \lambda_t,t)$. Our ansatz is a function of the form
\[h(x,\lambda, t) := \beta \exp(-\theta(r) t - \alpha(r) \lambda - r x).\]
We want this function to be in the domain of the generator $\mathcal{A}$ of the PDMP $(X, \lambda,\cdot)$. As already mentioned, the function $h$ has to be absolutely continuous, which is the case independent of the choice of the parameters $\alpha$, $r$ and $\theta$. Further it has to satisfy the integrability condition
\[\mathbb{E}_{(u,\lambda_0)} \left[ \sum_{i=1}^{\tilde N_t} \left\vert h(X_{T_i}, \lambda_{T_i}, T_i) - h(X_{T_i-}, \lambda_{T_i-}, T_i-)\right\vert \right] < \infty\] for all $t \geq 0$. Since $\tilde N$ is a Cox process with intensity $\lambda+\rho$, it is crucial to consider the whole sum and not only a single jump term. 
Additionally we want $h$ to satisfy the equation $\mathcal{A}h=0$ to ensure that it is a martingale. Obviously, this cannot hold for every choice of $r$, $\alpha$ and $\theta$. \\ \\
To motivate the explicit choice of our parameters let us assume that $h$ is in the domain and apply $\mathcal{A}$ to $h$. This gives us \begin{equation*}\begin{split}
\mathcal{A}h(x,&\lambda,t) = - \theta h(x,\lambda,t)-crh(x,\lambda,t) + \delta \lambda \alpha h(x,\lambda,t)  \\ &+\lambda h(x,\lambda,t) \int_0^\infty \! \left(e^{ru}-1\right) \, F_U(\mathrm{d}u) +\rho h(x,\lambda,t) \int_0^\infty \! \left(e^{-\alpha y}-1\right)\, F_Y(\mathrm{d}y)\overset{!}{=} 0.
\end{split}
\end{equation*}
Since $h$ is strictly positive we can reformulate the equation to \[\delta \lambda \alpha -cr-\theta +\lambda(M_U(r)-1) + \rho (M_Y(-\alpha)-1) = 0.\]  Here $M_U(s)$ and $M_Y(s)$ denote the moment-generating functions of the random variables $U$ and $Y$, which we assume to be finite.
The equation above has to hold for any $\lambda >0$, hence this is equivalent to \begin{align*}
&\delta \alpha + M_U(r)-1 =0, 
\\
& -cr-\theta + \rho (M_Y(-\alpha)-1) =0.
\end{align*}
Solving the above equations for some fixed $r$ we get the unique solutions $$\alpha(r) = \frac{1-M_U(r)}{\delta}$$ and $$\theta(r) = -cr + \rho \left(M_Y\left(\frac{M_U(r)-1}{\delta}\right)-1\right).$$
Now we still have to show that for this explicit choice of the parameters, the function $h$ is in the domain of the generator $\mathcal{A}$.

\begin{lemma}\label{lem_integrability}
Let $r$ be constant such that $M_U(r)$ is finite and define $\alpha(r):= \frac{1-M_U(r)}{\delta}$. Assume further that $M_Y(-\alpha(r))$ is finite. If $\theta (r) := -cr + \rho \left(M_Y(-\alpha(r))-1\right)$ and $\beta= \exp(ru+\alpha(r)\lambda_0)$, then $h(X_t,\lambda_t,t)$ is integrable and has expectation $1$ for all $t\geq 0$.
\begin{proof}
The expectation can be rewritten as
\begin{multline*}
\mathbb{E}_{(u,\lambda_0)}\left[ \beta \exp(-rX_t -\alpha(r)\lambda_t -\theta(r)t )\right] =
\exp(-rct -\theta(r) t + \alpha(r) \lambda_0)\\ \mathbb{E}_{(u,\lambda_0)}\left[\exp(-r \sum_{i=1}^{N_t} U_i - \alpha(r) \lambda_t)\right].
\end{multline*}
Conditioned on $\mathcal{F}^\lambda_t$, the counting process $N$ is an inhomogeneous Poisson process and as shown in \cite{Albrecher.2006} its integrated compensator has the form \[ \Lambda_t = \int_0^t \! \lambda_s \, \mathrm{d}s = \frac{1}{\delta}\left(\lambda_0 + \sum_{j=1}^{N^\lambda_t} Y_j - \lambda_t \right).\]
Using this we get
\begin{align*}
&\exp(-rct -\theta(r) t + \alpha(r) \lambda_0)\,\mathbb{E}_{(u,\lambda_0)}\left[\exp\left(-r \sum_{i=1}^{N_t} U_i - \alpha(r) \lambda_t\right)\right] = \\
&\exp(-rct -\theta(r) t + \alpha(r) \lambda_0)\,\mathbb{E}_{(u,\lambda_0)}\left[\exp\Big( \left( M_U(r)-1\right)\Lambda_t - \alpha(r) \lambda_t\Big)\right]
=\\ &\exp(-rct -\theta(r) t\,)\mathbb{E}\left[\exp\left(-\alpha(r) \sum_{j=1}^{N^\lambda_t} Y_j\right)\right].
\end{align*} 
The process $\sum_{j=1}^{N^\lambda_t} Y_j$ is a compound Poisson process, whose moment-generating function is $\exp(\rho t (M_Y(-\alpha(r))-1)$. By this and the definition of $\theta(r)$ we get that $h(X_t,\lambda_t,t)$ has expectation $1$.
\end{proof}
\end{lemma}
\begin{lemma}\label{lemma_is_in_domain}
Let all conditions of Lemma \ref{lem_integrability} be satisfied and assume further that \linebreak $\mathbb{E}\left[Y \exp\left(-\alpha(r) Y\right)\right]$ is finite. Then the function $h$ is in the domain of the generator.
\begin{proof}
Since $h$ is in $C^1$, we only have to check the integrability condition.
The jumps of $\lambda$ and $X$ have intensities, hence with probability $1$ they do not jump at the same time. Consequently, every jump time $T_i$ of $\tilde N$ can be uniquely identified with a jump time $T^\lambda_j$ of $N^\lambda$ or a jump time $T^X_k$ of $N$. Using this we get \begin{align*}
&\mathbb{E}_{(u,\lambda_0)} \left[ \sum_{i=1}^{\tilde N_t} \left\vert h(X_{T_i}, \lambda_{T_i}, T_i) - h(X_{T_i-}, \lambda_{T_i-}, T_i-)\right\vert \right] = \\ 
&\mathbb{E}_{(u,\lambda_0)} \left[ \sum_{i=1}^{ N_t} \left\vert h(X_{T^X_i}, \lambda_{T^X_i-}, T^X_i-) - h(X_{T^X_i-}, \lambda_{T^X_i-}, T^X_i-)\right\vert \right] 
\\& 
+\mathbb{E}_{(u,\lambda_0)} \left[ \sum_{i=1}^{ N^\lambda_t} \left\vert h(X_{T^\lambda_i-}, \lambda_{T^\lambda_i}, T^\lambda_i-) - h(X_{T^\lambda_i-}, \lambda_{T^\lambda_i-}, T^\lambda_i-)\right\vert \right].
\end{align*}
Let us now focus on the jumps in $\lambda$. 
\begin{align*}
&\mathbb{E}_{(u,\lambda_0)} \left[ \sum_{i=1}^{ N^\lambda_t} \left\vert h(X_{T^\lambda_i-}, \lambda_{T^\lambda_i}, T^\lambda_i-) - h(X_{T^\lambda_i-}, \lambda_{T^\lambda_i-}, T^\lambda_i-)\right\vert \right] = \\ 
&\mathbb{E}_{(u,\lambda_0)} \left[ \sum_{i=1}^{ N^\lambda_t}  h(X_{T^\lambda_i-}, \lambda_{T^\lambda_i-}, T^\lambda_i-)\left\vert \exp(-\alpha(r) Y_i)-1\right\vert \right]=\\&\mathbb{E}_{(u,\lambda_0)} \left[ \sum_{i=1}^{ N^\lambda_t}  h(X_{T^\lambda_i-}, \lambda_{T^\lambda_i-}, T^\lambda_i-)\mathbb{E} \left[\left\vert \exp(-\alpha(r) Y_i)-1\right\vert\, \Big| N^\lambda_t \cup\mathcal{F}_{T^\lambda_i-} \right]\right] .
\end{align*}
The triangle inequality and the fact that $Y_i$ is independent of $N^\lambda$ and $\mathcal{F}_{T^\lambda_i-}$ gives us
\begin{align*}
&\mathbb{E}_{(u,\lambda_0)} \left[ \sum_{i=1}^{ N^\lambda_t}  h(X_{T^\lambda_i-}, \lambda_{T^\lambda_i-}, T^\lambda_i-)\left\vert \exp(-\alpha(r) Y_i)-1\right\vert \right] \leq \\
&\mathbb{E}_{(u,\lambda_0)} \left[ \sum_{i=1}^{ N^\lambda_t}  h(X_{T^\lambda_i-}, \lambda_{T^\lambda_i-}, T^\lambda_i-)(M_Y(-\alpha(r))+1) \right] =\\
&(M_Y(-\alpha(r))+1)\,\mathbb{E}_{(u,\lambda_0)} \left[\int_{(0,t]} \! h(X_{s-}, \lambda_{s-}, s-) \, \mathrm{d}N^\lambda_s\right].
\end{align*}
Compensating with the intensity $\lambda$ yields\begin{align*}
 &(M_Y(-\alpha(r))+1)\,\mathbb{E}_{(u,\lambda_0)} \left[\int_{(0,t]} \! h(X_{s-}, \lambda_{s-}, s-) \, \mathrm{d}N^\lambda_s\right] = \\
&(M_Y(-\alpha(r))+1)\,\mathbb{E}_{(u,\lambda_0)} \left[\int_0^t \rho h(X_{s-},\lambda_{s-},s-) \, \mathrm{d}s\right]=\\
&(M_Y(-\alpha(r))+1)\,\mathbb{E}_{(u,\lambda_0)} \left[\int_0^t \rho h(X_s,\lambda_s,s) \, \mathrm{d}s\right].
\end{align*}
Since $h$ is positive, we can apply Tonelli's theorem to interchange expectation and integration and get the result $(M_Y(-\alpha(r))+1)\rho t$ which is finite for all $t \geq 0$.
\\ 
Considering the jumps in $X$, the same arguments as before result in \begin{align*}
&\mathbb{E}_{(u,\lambda_0)} \left[ \sum_{i=1}^{ N_t} \left\vert h(X_{T^X_i}, \lambda_{T^X_i-}, T^X_i-) - h(X_{T^X_i-}, \lambda_{T^X_i-}, T^X_i-)\right\vert \right] \leq \\
&\quad(M_U(r)+1) \int_0^t \! \mathbb{E}_{(u,\lambda_0)} \left[\lambda_s h(X_s,\lambda_s,s)\right] \, \mathrm{d}s.
\end{align*}\\
Now take a look at $\mathbb{E}_{(u,\lambda_0)} \left[\lambda_t h(X_t,\lambda_t,t) \right]$. We can use the same ideas as in the prove of the integrability of $h(X_t,\lambda_t,t)$ to obtain
\[ \mathbb{E}_{(u,\lambda_0)} \left[\lambda_t h(X_t,\lambda_t,t) \right] = \exp(-rct -\theta(r) t)\mathbb{E}_{(u,\lambda_0)}\left[\lambda_t\exp\left(-\alpha(r) \sum_{j=1}^{N^\lambda_t} Y_j\right)\right]. \]
The drift of $\lambda$ is only decreasing, hence $\lambda_t \leq \sum_{i=1}^{N^\lambda_t} Y_i$. Using this and the independence of $N^ \lambda$ and the $Y_i$ we get by conditioning on $N^\lambda_t$ \[\mathbb{E}_{(u,\lambda_0)} \left[\lambda_t h(X_t,\lambda_t,t) \right] \leq \mathbb{E}\left[\sum_{i=1}^{N^\lambda_t} \prod_{j=1}^{N^\lambda_t} \mathbb{E}\left[Y_i \exp(-\alpha(r) Y_j)\right]\right]. \]
For $i \neq j$ the random variables $Y_i$ and $Y_j$ are independent hence $\mathbb{E} \left[Y_i \exp(-\alpha(r) Y_j)\right] = \mathbb{E} \left[Y\right]M_Y(-\alpha(r))$. For $i=j$ we assumed that $\mathbb{E}\left[Y \exp(-\alpha(r)Y)\right]$ is finite.
By this we have that \begin{align*}
\mathbb{E}_{(u,\lambda_0)} \left[\lambda_t h(X_t,\lambda_t,t) \right] \leq \frac{\mathbb{E} \left[Y \exp(-\alpha(r)Y)\right]}{\mathbb{E}\left[Y \right]M_Y(-\alpha(r))} \mathbb{E} \left[N^\lambda_t \exp(N^\lambda_t \log\left(\mathbb{E} \left[Y \right]M_Y(-\alpha(r))\right) \right].
\end{align*}
The random variable $N^\lambda_t$ is Poisson distributed with parameter $\rho t$ hence \begin{align*}\mathbb{E}\left[N^\lambda_t \exp(N^\lambda_t \log\left(\mathbb{E}\left[Y \right]M_Y(-\alpha(r))\right) \right] = \rho t &\mathbb{E}\left[Y \right]M_Y(-\alpha(r))\\ &\exp\left(\rho t\left(\mathbb{E}\left[Y \right]M_Y(-\alpha(r))-1\right)\right).
\end{align*}
Using the derived results we know that there are some positive constants $K$ and $\gamma$ such that 
\begin{align*}
&\mathbb{E}_{(u,\lambda_0)} \left[ \sum_{i=1}^{ N_t} \left\vert h(X_{T^X_i}, \lambda_{T^X_i-}, T^X_i-) - h(X_{T^X_i-}, \lambda_{T^X_i-}, T^X_i-)\right\vert \right]  \\ \\
&\leq \int_0^t K\exp(\gamma s) s \, \mathrm{d}s = K \frac{\exp(\gamma t)(\gamma t -1)+1}{\gamma^2} < \infty.
\end{align*}
Consequently, the function $h$ is in the domain of the generator.
\end{proof}
\end{lemma}
These results lead immediately to the following Theorem.
\begin{theorem}
Under the assumptions of Lemma \ref{lemma_is_in_domain}, the process $M^r_t:=h(X_t,\lambda_t,t)$ is an expectation $1$ martingale.
\begin{proof}
By Lemma \ref{lem_integrability}, the process is integrable and has constant expectation $1$. The function $h$ is in the domain of the generator $\mathcal{A}$ and satisfies $\mathcal{A}h =0$. Therefore $M^r_t$ is a martingale.
\end{proof}
\end{theorem}
Using these martingales we can define a family of measures $\mathbb{Q}^{(r)}$ such that \[ \left.\frac{\mathrm{d\mathbb{Q}^{(r)}}}{\mathrm{d}\mathbb{P}_{(u,\lambda_0)}}\right\vert_{\mathcal{F}_t} = M^{(r)}_t.\] The exponential form of the change of measure allows us to exploit of the results shown in \cite{ZbigniewPalmowski.2002} to derive the behaviour of the combined process under the new measures $\mathbb{Q}^{(r)}$.
\begin{lemma}
Let $r$ be such that $M^r$ is well defined. Then under the measure $\mathbb{Q}^{(r)}$, the process $(X, \lambda, t
\cdot)$ is again a PDMP with generator \begin{align*}
\mathcal{A}^{(r)} f(x,\lambda,t) = &c f_x(x,\lambda,t) - \delta \lambda f_\lambda(x,\lambda,t) + f_t(x,\lambda,t) \\&+ \lambda \int_0^\infty \! (f(x-u,\lambda,t) -f(x,\lambda,t)) e^{ru} \, F_U(\mathrm{d}u) \\&+ \rho \int_0^\infty \! (f(x, \lambda+y, t) - f(x,\lambda,t) ) e^{-\alpha(r) y} \, F_Y(\mathrm{d}y),
\end{align*} 
\end{lemma} 
So far, we have found a new family of measures but we have to identify a measure that fits our needs. Motivated by the definition of the adjustment coefficient in the classical model we consider the function $\theta(r) = -cr + \rho (M_Y(-\alpha(r)) -1 )$.

\begin{lemma}
The function $\theta(r)$ is convex and satisfies $\theta(0)=0$. 
\begin{proof}
To show convexity we use the fact, that moment-generating functions are log-convex and therefore convex. Moreover they are twice differentiable. Consequently $\theta$ is twice differentiable too and its derivatives are \begin{align*}
\theta'(r)& = -c + \frac{\rho}{\delta} M_Y'(\frac{M_U(r)-1}{\delta}) M_U'(r) ,\\ \\
\theta''(r) &= \frac{\rho}{\delta^2} M_Y'\left(\frac{M_U(r)-1}{\delta}\right) M_U'(r)^2 + \frac{\rho}{\delta} M_Y'\left(\frac{M_U(r)-1}{\delta}\right)M_U''(r).
\end{align*}
By convexity of the moment-generating functions, we know that their second derivatives are non-negative. To ensure that $\theta$ is convex, we have to check if the first derivative of the MGF of $Y$ is non-negative too. Equivalently we show that the MGF of $Y$ is monotone increasing. Let now $r>s$ then $ \mathbb{E}\left[e^{rY}\right] = \mathbb{E}\left[e^{sY} e^{(r-s)Y}\right]. $ The random variable $Y$ is almost surely positive and $r-s$ is positive too. Hence $e^{(r-s)Y} > 1$ almost surely. This gives us \[ M_Y(r) = 
\mathbb{E}\left[e^{sY} e^{(r-s)Y}\right] >  \mathbb{E}\left[e^{sY}\right] =M_Y(s).\] Consequently the first derivative of $M_Y(r)$ is non-negative. Therefore $\theta$ is convex and since $M_U(0)=M_Y(0)=1$ we get that $\theta(0)=0$.
\end{proof}
\end{lemma}

\begin{lemma}

Let $r$ be such that the measure $\mathbb{Q}^{(r)}$ is well defined and assume there is some $\varepsilon >0$ such that $M_U(r+\varepsilon)$ and $M_Y(-\alpha(r) + \varepsilon)$ are finite. Then \[ \lim_{t \to \infty} \, \frac{\mathbb{E}^{\mathbb{Q}^{(r)}}\left[X_t\right]}{t} = - \theta'(r). \]
\begin{proof}
To show this property, we can use the ideas of the proof of Lemma \ref{lemma_expectation_surplus}. The main difference is, that we apply the generator $\mathcal{A}^{(r)}$. Again we obtain \[\mathbb{E}^{\mathbb{Q}^{(r)}}\left[X_t\right]=u+ ct -M_U(r) \mathbb{E}^{\mathbb{Q}^{(r)}}\left[U\right] \int_0^t \! \mathbb{E}^{\mathbb{Q}^{(r)}}\left[\lambda_s\right] \, \mathrm{d}s.\]
The expectation of $\lambda_t$ under $\mathbb{Q}^{(r)}$ satisfies \[\mathbb{E}^{\mathbb{Q}^{(r)}}\left[\lambda_t\right] = \frac{\rho}{\delta} M_Y(-\alpha(r)) \mathbb{E}^{\mathbb{Q}^{(r)}}\left[Y\right] (1- e^{-\delta t}) + e^{-\delta t} \lambda_0.\]

The expectations $ \mathbb{E}^{\mathbb{Q}^{(r)}}\left[U\right]$ and $\mathbb{E}^{\mathbb{Q}^{(r)}}\left[Y\right]$ can easily be obtained from \[ M^{\mathbb{Q}^{(r)}}_U(s) = \frac{M_U(s+r)}{M_U(r)} \] and \[ M^{\mathbb{Q}^{(r)}}_Y(s) = \frac{M_Y(s-\alpha(r))}{M_Y(-\alpha(r))}. \] Consequently, $
\mathbb{E}^{\mathbb{Q}^{(r)}}\left[U\right] = \frac{M_U'(r)}{M_U(r)} $ and $
\mathbb{E}^{\mathbb{Q}^{(r)}}\left[Y\right] =  \frac{M_Y'(-\alpha(r))}{M_Y(-\alpha(r))}$. \\
Combining these results we get 
\begin{align*}
\lim_{t \to \infty} \, \frac{\mathbb{E}^{\mathbb{Q}^{(r)}}\left[X_t\right]}{t} &= c-\frac{\rho}{\delta}M_Y'(-\alpha(r))M_U'(r) = - \theta'(r).
\end{align*}
\end{proof}
\end{lemma}

\begin{assumption}\label{assumption_epsilon}
From now on we assume that there exists a positive solution $R$ to the equation $\theta(R) =0$, that $\mathbb{Q}^{(R)}$ is well defined and that for some $\varepsilon>0$ both $M_U(R+\varepsilon)$ and $M_Y(\varepsilon-\alpha(R))$ are finite.
\end{assumption}
This assumption ensures that the measure $\mathbb{Q}^{(R)}$ is well defined and that we can express the expectation of $Y$ and $U$ in terms of their original moment-generating functions. One example where this is satisfied is the following.
\begin{example}\label{example_1}
Let $\mu$ and $\kappa$ be positive constants. If $Y \sim Exp(\mu)$ and $U \sim Exp(\kappa)$, the net profit condition simplifies to $c > \frac{\rho}{\delta \kappa \mu}$. The moment-generating functions are given by $M_U(r)= \frac{\kappa}{\kappa-r} $ and $M_Y(-\alpha(r))= \frac{\mu}{\mu+\alpha(r)},$ where $r<\kappa$ and $-\alpha(r) < \mu$.
If we fix some $r < \frac{\mu \delta \kappa}{1+\delta\mu}$ we can determine the functions $ \alpha(r) = - \frac{r}{\delta (\kappa -r)}$ and \[\theta(r) = -cr + \rho\left(\frac{r}{\mu\delta(\kappa-r)-r}\right).\]
Solving the equation $\theta(r) = 0$ gives us the solutions $r_1=0$ and \[ R:=r_2 = \frac{\mu\delta\kappa c - \rho}{(1+\mu \delta)c},\] which is positive by the net profit condition. Now we want to show that there is some $\varepsilon>0$ such that $R+\varepsilon < \frac{\mu\delta}{1+\mu\delta}\kappa$ and $\varepsilon-\alpha(R) < \mu$. The first inequality is equivalent to \[\varepsilon < \frac{\rho}{(1+\mu \delta)c},\] which is a strictly positive upper bound. The second condition can be rewritten to \[ \varepsilon < \frac{\mu\rho \delta + \rho}{\delta\kappa c + \rho \delta},\] which is positive too. Consequently, Assumption \ref{assumption_epsilon} is satisfied.
\end{example}
\begin{lemma}
For every $u\geq0$ and $\lambda_0>0$ we have that $\mathbb{Q}^{(R)} \left[ \tau_u < \infty\right] =1.$
\begin{proof}
We already know that $\lim_{t \to \infty} \, \frac{\mathbb{E}^{\mathbb{Q}^{(R)}}\left[X_t\right]}{t} = - \theta'(R) $ holds true. If we can show that $\theta'(R) >0$ then ruin occurs almost surely under the new measure.
The function $\theta$ is convex and satisfies $\theta(0)=\theta(R)=0$. Further we have that \[\theta'(0)= -c + \frac{\rho}{\delta} \mathbb{E}\left[Y\right]\mathbb{E}\left[U\right], \] which is smaller than $0$ by the net profit condition. Therefore there exists $0<r<R$ such that $\theta(r) <0$. Since $\theta(R)>\theta(r)$ it follows by the mean-value theorem that there is a $\tilde r \in (r,R)$ such that \[ \theta'(\tilde r) = \frac{\theta(R)-\theta(r)}{R-r} >0.\] By convexity we know that $\theta'$ is a monotone increasing function and $\theta'(R)\geq \theta'(\tilde r)>0$.

\end{proof}
\end{lemma}
Similar to the classical model and the Björk-Grandell model which is considered in \cite{Schmidli.1997}, we have found a new measure under which ruin occurs almost surely. We can use this to get an upper bound for the ruin probability.
\begin{theorem}
Under our assumptions \[ \psi(u,\lambda_0) \leq e^{-\alpha(R)\lambda_0} e^{-Ru}.\]
\begin{proof}
The ruin probability can be rewritten as \begin{align*}
\psi(u,\lambda_0) &= \mathbb{E}_{(u,\lambda_0)}\left[I_{\left\lbrace \tau_u < \infty\right\rbrace} \right]= \mathbb{E}^{\mathbb{Q}^{(R)}}\left[I_{\left\lbrace \tau_u < \infty\right\rbrace} \left(M_{\tau_u}^R\right)^{-1}\right] \\ &= \exp(-Ru-\alpha(R) \lambda_0) \mathbb{E}^{\mathbb{Q}^{(R)}}\left[\exp\left(RX_{\tau_u} + \alpha(R) \lambda_{\tau_u}\right) \right].
\end{align*}
By definition of $\tau_u$, the value $X_{\tau_u}$ is negative and since $R>0$ we have that $M_U(R)>1$. Consequently, $\alpha(R) <0$. By this we get that $\exp(RX_{\tau_u} + \alpha(R) \lambda_{\tau_u}) \leq 1$ and \[ \psi(u,\lambda_0) \leq \exp(-Ru-\alpha(R) \lambda_0).\]
\end{proof}
\end{theorem}
\section{The Renewal Equation}\label{section_renewal_equation}
We now want to use Theorem \ref{renewal_theorem_schmidli} to get information about the asymptotic behaviour of the ruin probability $\psi(u,\lambda_0)$. Because of the dependence on $\lambda_0$, we have to choose the renewal times $\left\lbrace S_+(i)\right\rbrace_{i \in \mathbb{N}}$ such that $\lambda_{S_+(i)}=\lambda_0$. To exploit the renewal equation we have to ensure that there are infinitely many renewal times and that they are almost surely finite. For this, we will use the ideas from \cite{Orsingher.1982} to get an intensity for the number of upcrossings of the process $\lambda$ through some level $l$.
\begin{lemma}
Let $\lambda$ be the Markovian shot-noise process and $l>0$ arbitrary. The process counting all upcrossings of $\lambda$ through $l$ has intensity \[\nu^+_l(t) = \rho \int_0^l (1-F_Y(l-z)) \, F_\lambda(\mathrm{d}z,t),\] where $F(z,t)= \mathbb{P}_{(u,\lambda_0)}\left[\lambda_t \leq z\right]$ is the CDF of $\lambda_t$.

\begin{proof}
Consider for some small $\Delta t$ the probability $\mathbb{P}_{(u,\lambda_0)} \left[ \lambda_t \leq l, \lambda_{t+\Delta t} > l \right]$. The jumps of $\lambda$ are governed by a Poisson process with rate $\rho$ hence \begin{align*}
&\mathbb{P}_{(u,\lambda_0)} \left[ \lambda_t \leq l, \lambda_{t+\Delta t} > l \right] = \mathbb{P}_{(u,\lambda_0)} \left[ N^{\lambda}_{t +\Delta t} - N^\lambda_t = 1, 
\lambda_t \leq l, \lambda_{t+\Delta t} > l \right] + o(\Delta t)= \\ 
& \mathbb{P}_{(u,\lambda_0)} \left[ N^{\lambda}_{t +\Delta t} - N^\lambda_t = 1,  \lambda_t \leq l, \lambda_t e^{-\delta \Delta t} + Y e^{-\delta (t+\Delta t -T)} > l \right] + o(\Delta t)= \\ 
& \mathbb{P}_{(u,\lambda_0)} \left[ N^{\lambda}_{t +\Delta t} - N^\lambda_t = 1,  \lambda_t \leq l, \lambda_t + Y e^{-\delta (t -T)} > le^{\delta \Delta t} \right] + o(\Delta t).
\end{align*}
Here $T$ denotes the jump time occurring between $t$ and $t+\Delta t$ and $Y$ is the corresponding shock. The random time $T-t$ can be represented as $\Theta \Delta t$, where $\Theta$ is a random variable which takes values in the interval $(0,1)$. Consequently we have that  $Ye^{\delta \Theta \Delta t} \in (Y,Ye^{\delta \Delta t})$. Using this we can bound the above probability by \begin{align*}
&\mathbb{P}_{(u,\lambda_0)} \left[ N^{\lambda}_{t +\Delta t} - N^\lambda_t = 1,  \lambda_t \leq l, \lambda_t + Y e^{\delta \Delta t} > le^{\delta \Delta t} \right] + o(\Delta t) \geq\\ &\mathbb{P}_{(u,\lambda_0)} \left[ N^{\lambda}_{t +\Delta t} - N^\lambda_t = 1,  \lambda_t \leq l, \lambda_t + Y e^{-\delta \Theta \Delta t} > le^{\delta \Delta t} \right] + o(\Delta t) \geq \\ &\mathbb{P}_{(u,\lambda_0)} \left[ N^{\lambda}_{t +\Delta t} - N^\lambda_t = 1,  \lambda_t \leq l, \lambda_t + Y > le^{\delta \Delta t} \right] + o(\Delta t).
\end{align*}

Let us focus on the upper bound. The term $N^\lambda_{t+\Delta t} - N^\lambda_t$ is independent of $\lambda_t$ and $Y$ so we can rewrite \begin{align*}
&\mathbb{P}_{(u,\lambda_0)} \left[ N^{\lambda}_{t +\Delta t} - N^\lambda_t = 1,  \lambda_t \leq l, \lambda_t + Y e^{\delta \Delta t} > le^{\delta \Delta t} \right] + o(\Delta t) = \\&\rho \Delta t \mathbb{P}_{(u,\lambda_0)} \left[ \lambda_t \leq l, \lambda_t + Y e^{\delta \Delta t} > le^{\delta \Delta t} \right] + o(\Delta t) = \\ & \rho \Delta t\, \mathbb{E}_{(u,\lambda_0)} \left[\mathbb{E}_{(u,\lambda_0)} \left[ I_{\left\lbrace\lambda_t \leq l\right\rbrace} I_{\left\lbrace Y> l- \lambda_t e^{-\delta \Delta t}\right\rbrace}\,\left\vert \, \lambda_t \right.\right]\right] + o(\Delta t) = \\ & \rho \Delta t \int_0^l \! (1-F_Y(l-ze^{-\delta \Delta t})) \, F_\lambda(\mathrm{d}z, t) + o(\Delta t).
\end{align*}
Now, let us divide by $\Delta t$ and consider the limit of $\Delta t \to 0$. Since $F_Y(l-ze^{-\delta \Delta t})$ decreases as $\Delta t$ becomes smaller, we get by the right continuity of CDFs that this tends to \[\rho \int_0^l \! (1-F_Y(l-z)) \, F_\lambda(\mathrm{d}z, t).\] Using the same arguments we can show that the lower bound divided by $\Delta t$ converges to the same value. Hence, the term $\frac{1}{\Delta t} \mathbb{P}_{(u,\lambda_0)} \left[ \lambda_t \leq l, \lambda_{t+\Delta t} > l \right]$ converges too.
\end{proof}
\end{lemma}

\begin{assumption} \label{assumption_crossings}
From now on we assume that \[  \int_0^\infty \int_0^{\lambda_0} (1-F_Y^{\mathbb{Q}^{(R)}}(\lambda_0-z)) \,F^{\mathbb{Q}^{(R)}}_\lambda(\mathrm{d} z,t) \, \mathrm{d}t = \infty, \] where $F^{\mathbb{Q}^{(R)}}_\lambda(z,t):= \mathbb{Q}^{(R)}\left[\lambda_t \leq z\right]$ and $F_Y^{\mathbb{Q}^{(R)}}(x) = \mathbb{Q}^{(R)}\left[Y \leq x\right]$.
\end{assumption}

This assumption guarantees that there are infinitely many upcrossings of the process through $\lambda_0$ under the measure $\mathbb{Q}^{(R)}$. The structure of our Markovian shot-noise process gives us, that upcrossings can only happen through shock events and downcrossings are due to the continuous drift. Consequently, there have to be infinitely many continuous downcrossings and recurrence times $\left\lbrace S(i)\right\rbrace_{i \in \mathbb{N}},$ such that $\lambda_{S(i)} = \lambda_0$.\\
One example which satisfies Assumption \ref{assumption_crossings} is the following.
\begin{example}
Consider the same configuration as in Example \ref{example_1}. Under the new measure $\mathbb{Q}^{(R)}$, the shocks are again exponentially distributed with parameter $\mu+\alpha(R)$ and the new jumping intensity of $\lambda$ is \[\tilde{\rho} =\rho M_Y(-\alpha(R))= \frac{\mu\delta\kappa c+ \mu \delta \rho}{\mu \delta +1}.\] 
Assume that $\frac{\tilde \rho}{\delta}=n \in \mathbb{N}.$ Like in \cite{Orsingher.1982} we can determine the distribution of $Y(t)$ using its characteristic function \[ K_t(s) = \mathbb{E}^{\mathbb{Q}^{(R)}}\left[\exp(is \lambda(t)) \right] = \left(e^{-\delta t} +(1-e^{\delta t})\frac{\mu+\alpha(R)}{(\mu+\alpha(R))-is}\right)^n.\]
This is the characteristic function of the random variable $\eta = \sum_{i=1}^{B_t} Y_i$, where $$B_t \sim \mathit{B}\left(n,1-e^{-\delta t}\right).$$
Consequently, $\lambda(t)$ admits a density of the form \[f(z,t)= \sum_{j=1}^n \binom{n}{j} e^{-\delta t (n-j)} (1-e^{-\delta t})^j (\mu+\alpha(R))^j e^{-(\mu+\alpha(R)) z}\frac{z^{j-1}}{(j-1)!}.\]
Using this, the intensity of the upcrossings is given by \[ \nu^+_{\lambda_0}(t) = \rho \sum_{j=0}^n \binom{n}{j}e^{-\delta t (n-j)} (1-e^{-\delta t})^j \frac{(\mu+\alpha(R))^j \lambda_0^j}{j!}e^{-(\mu+\alpha(R))\lambda_0}. \] Since $\frac{(\mu+\alpha(R))^j \lambda_0^j}{j!}$ has a positive lower bound $\tilde c$ we get that \begin{align*}
\int_0^\infty \nu^+_{\lambda_0} (t) \, \mathrm{d}t \geq  \int_0^\infty \rho \tilde{c}e^{-(\mu+\alpha(R))\lambda_0} \mathrm{d}t = \infty.
\end{align*}
\end{example} \vspace{1cm}
If Assumption \ref{assumption_crossings} holds we have that under the measure $\mathbb{Q}^{(R)}$, the surplus process tends to $-\infty$ and $\lambda$ returns to $\lambda_0$ infinitely often. Hence, we can define a sequence of renewal times $\left\lbrace S_+(i) \right \rbrace_{i \in \mathbb{N}_0}$ via $S_+(0)=0$ and  $S_+(i) = \min \left\lbrace S(i)>S_+(i-1) \, \left\vert \, X_{S(i)} < X_{S_+(i-1)} \right.\right\rbrace$ which satisfies $\mathbb{Q}^{(R)} \left[ S_+(i) < \infty\right] =1$ for all $i$. We will use these renewal times similar to the ladder epochs in the classical ruin model.\\
Define $$B(x) = \mathbb{P}_{(u,\lambda_0)}\left[S_+(1) < \infty,\, u-X_{S_+(1)} \leq x \right]$$ and $$p(u,x) = \mathbb{P}_{(u,\lambda_0)}\left[ \tau_u \leq S_+(1) \, | \, S_+(1) < \infty,\, X_{S_+(1)}=u-x \right].$$ Then the ruin probability satisfies:
\[ \psi(u,\lambda_0) = \int_0^u \psi(u-x,\lambda_0) (1-p(u,x)) \, B(\mathrm{d}x) + \mathbb{P}_{(u,\lambda_0)}\left[\tau_u \leq S_+(1),\, \tau_u < \infty\right]. \] 
This may look like a renewal equation but the distribution $B$ is defective. We solve this problem by multiplying both sides with $e^{Ru}$, which is equivalent to a measure change from $\mathbb{P}$ to $\mathbb{Q}^{(R)}$, and obtain:
\begin{multline}\label{equation_renewal_our_case}
 \psi(u,\lambda_0)e^{Ru} = \int_0^u \psi(u-x,\lambda_0)e^{R(u-x)} (1-p(u,x))e^{Rx} \, B(\mathrm{d}x) \\+ \mathbb{P}_{(u,\lambda_0)}\left[\tau_u \leq S_+(1),\, \tau_u < \infty\right]e^{Ru}.
\end{multline} 
\begin{lemma}
The distribution $\tilde B$ defined by $\tilde B (\mathrm{d}x) = e^{Rx} B(\mathrm{d}x)$ is non-defective.
\begin{proof}
Using the definition of $\tilde B$ we get \[
\int_\mathbb{R}\tilde B(\mathrm{d}x) = \int_\mathbb{R} e^{Rx} B(\mathrm{d}x) = \mathbb{E}_{(u,\lambda_0)} \left[ e^{R(u-X_{S_+(1)}}I_{\left\lbrace S_+(1) < \infty\right\rbrace} \right]. \]
Now focus on our martingale $M^R$ at time $S_+(1)$ and observe that \begin{align*}
M_{S_+(1)}^R = \exp\left(\alpha(R) \lambda_0 + Ru-\alpha(R) \lambda_{S_+(1)} - RX_{S_+(1)}\right) = \exp\left(R(u-X_{S_+(1)})\right)
\end{align*}
Using this leads to \[\int_\mathbb{R} \tilde B(\mathrm{d}x) = 
\mathbb{E}^{\mathbb{Q}^{(R)}}\left[I_{\left\lbrace S_+(1) < \infty \right\rbrace}\right] = \mathbb{Q}^{(R)}\left[S_+(1) < \infty \right] =1.\] Consequently, $\tilde B$ is not defective.
\end{proof}
\end{lemma}
Even though we have found a renewal equation, we still have to show that all functions appearing in Equation \ref{equation_renewal_our_case} satisfy the assumptions of Theorem \ref{renewal_theorem_schmidli}. 
\begin{assumption}\label{assumption_r}
From now on we assume that there exists an $\varepsilon>0$ such that for $r:= (1+\varepsilon)R$ the measure $\mathbb{Q}^{(r)}$ is well defined and $$\mathbb{E}_{(u,\lambda_0)}\left[ e^{-r(X_{S_+(1)}-u)} I_{\left\lbrace S_+(1) < \infty\right\rbrace}\right] < \infty.$$
\end{assumption}
Since $S_+(1)$ depends on $X$ and $\lambda$, this assumption may be hard to check. Alternatively, we can use the following lemma, which allows us to focus on the first recurrence time $S(1)$.
\begin{lemma}
Let $\varepsilon>0$ such that for $r:=(1+\varepsilon)R$ the measure $\mathbb{Q}^(r)$ is well defined. Then \[ \mathbb{E}_{(u,\lambda_0)}\left[\exp(-r(X_{S_+(1)}-u))I_{\left\lbrace S_+(1) < \infty \right\rbrace}\right] <\infty \] if and only if \[ \mathbb{E}_{(u,\lambda_0)}\left[\exp\left(-r(X_{S(1)}-u)\right)I_{\left\lbrace S(1)<\infty\right\rbrace}\right] < \infty.\]
\begin{proof}
At first assume that \[ \mathbb{E}_{(u,\lambda_0)}\left[\exp(-r(X_{S_+(1)}-u))\right] <\infty \] holds. By definition $S_+(1) \geq S(1)$ and $\theta(r)>0$. Consequently \begin{align*}
\mathbb{E}_{(u,\lambda_0)}\left[\exp\left(-r(X_{S(1)}-u)\right)I_{\left\lbrace S(1)<\infty\right\rbrace}\right] &= \mathbb{E}^{\mathbb{Q}^{(r)}}\left[\exp(\theta(r)S(1))\right] \leq \mathbb{E}^{\mathbb{Q}^{(r)}}\left[\exp(\theta(r)S_+(1))\right] \\&=  \mathbb{E}_{(u,\lambda_0)}\left[\exp(-r(X_{S_+(1)}-u))I_{\left\lbrace S_+(1) < \infty \right\rbrace}\right] <\infty
\end{align*}
Let us now assume that \[ \mathbb{E}_{(u,\lambda_0)}\left[\exp\left(-r(X_{S(1)}-u)\right)I_{\left\lbrace S(1)<\infty \right\rbrace }\right]=:C < \infty\] holds true. Then \begin{align*}
\mathbb{E}_{(u,\lambda_0)}\left[\exp(-r(X_{S_+(1)}-u))I_{\left\lbrace S_+(1) < \infty\right\rbrace}\right] &= \mathbb{E}^{\mathbb{Q}^{(R)}}\left[\exp(-\varepsilon R(X_{S_+(1)}-u))\right] \\&= \sum_{i=1}^\infty \mathbb{E}^{\mathbb{Q}^{(R)}}\left[\exp(-\varepsilon R(X_{S(i)}-u))I_{\left \lbrace S_+(1)=S(i)\right\rbrace}\right].
\end{align*}
The indicator can be split up to \[ I_{\left\lbrace S_+(1)=S(i)\right\rbrace} =I_{\left\lbrace S_+(1) >S(i-1)\right\rbrace}I_{\left\lbrace X_{S(i)} < u\right\rbrace}= \prod_{j=1}^{i-1} I_{\left\lbrace X_{S(j)} \geq u\right\rbrace} I_{\left\lbrace X_{S(i)} < u\right\rbrace}.\] Observe further that, with $S(0)=0$, the random variables $$(\xi_j)_{j \geq 1}:=(X_{S(j)}-X_{S(j-1)})_{j \geq 1}$$ are i.i.d.. Consequently, $X_{i-1}-u= \sum_{j=1}^{i-1} \xi_j$ holds true for all $i$.
Using this, we get \begin{align*}
\mathbb{E}^{\mathbb{Q}^{(R)}}\left[\exp(-\varepsilon R(X_{S(i)}-u))I_{\left \lbrace S_+(1)=S(i)\right\rbrace}\right] \leq \mathbb{E}^{\mathbb{Q}^{(R)}}&\left[\exp\left(-\varepsilon R\sum_{j=1}^{i-1} \xi_j\right)I_{\left\lbrace \sum_{j=1}^{i-1} \xi_j >0\right\rbrace}\right. \\ 
&\left. \mathbb{E}^{\mathbb{Q}^{(R)}}\left[\left. \exp\left( -\varepsilon R \xi_i\right) I_{\left\lbrace X_{S(i)} < u \right\rbrace }\right\vert\sum_{j=1}^{i-1} \xi_j\right]\right]
\end{align*}
Let us focus on the conditional expectation. The indicator is less or equal to $1$ and $\xi_i$ is independent of the condition. Hence \[\mathbb{E}^{\mathbb{Q}^{(R)}}\left[ \left.\exp\left( -\varepsilon R \xi_i\right) I_{\left\lbrace X_{S(i)} < u \right\rbrace}\right\vert\sum_{j=1}^{i-1} \xi_j\right] \leq \mathbb{E}^{\mathbb{Q}^{(R)}}\left[ \exp\left( -\varepsilon R \xi_i\right)\right] = C< \infty.\]
By this we get that \[\mathbb{E}^{\mathbb{Q}^{(R)}}\left[\exp(-\varepsilon R(X_{S(i)}-u))I_{\left \lbrace S_+(1)=S(i)\right\rbrace}\right] \leq C\, \mathbb{E}^{\mathbb{Q}^{(R)}}\left[\exp\left(-\varepsilon R\sum_{j=1}^{i-1} \xi_j\right)I_{\left\lbrace \sum_{j=1}^{i-1} \xi_j >0\right\rbrace}\right].\]
Now we want to bound the remaining expectation. For this we observe that for all $\tilde \varepsilon >0$ \[ \mathbb{E}^{\mathbb{Q}^{(R)}}\left[\exp\left(-\varepsilon R\sum_{j=1}^{i-1} \xi_j\right)I_{\left\lbrace \sum_{j=1}^{i-1} \xi_j >0\right\rbrace}\right] \leq \mathbb{E}^{\mathbb{Q}^{(R)}}\left[\exp\left(\tilde\varepsilon R\sum_{j=1}^{i-1} \xi_j\right)\right].\]
To choose $\tilde \varepsilon$ in a suitable way, we focus on the properties of $\theta$. This function is convex and satisfies $\theta(0)=\theta(R) =0$ and $\theta'(0)<0$. Consequently, there exists a $\tilde r \in (0,R)$ such that $\theta(\tilde r) <0$ Choosing $\tilde \varepsilon = 1-\frac{\tilde r}{R} \in (0,1)$ we have that \begin{align*}
&\mathbb{E}^{\mathbb{Q}^{(R)}}\left[\exp\left(-\varepsilon R\sum_{j=1}^{i-1} \xi_j\right)I_{\left\lbrace \sum_{j=1}^{i-1} \xi_j >0\right\rbrace}\right] \leq \mathbb{E}^{\mathbb{Q}^{(R)}}\left[\exp\left(\tilde\varepsilon R\sum_{j=1}^{i-1} \xi_j\right)\right] = \mathbb{E}^{\mathbb{Q}^{(R)}}\left[\exp\left(\tilde\varepsilon R\xi_1\right)\right]^{i-1} 
\\[0.3cm]& = \mathbb{E}_{(u,\lambda_0)}\left[\exp\left((\tilde\varepsilon-1) R(X_{S(1)}-u)\right)I_{\left\lbrace S(1) < \infty\right\rbrace}\right]^{i-1}= \mathbb{E}_{(u,\lambda_0)}\left[\exp\left(-\tilde r (X_{S(1)}-u)\right)I_{\left\lbrace S(1) < \infty\right\rbrace}\right]^{i-1}
\\[0.3cm]& = \mathbb{E}^{\mathbb{Q}^{(\tilde r)}}\left[\exp\left(\theta(\tilde r)S(1)\right)I_{\left\lbrace S(1) < \infty\right\rbrace}\right]^{i-1}.
\end{align*} 
By construction we have that $\theta(\tilde r) < 0$ and $S(1)>0$ hence $$\mathbb{E}^{\mathbb{Q}^{(\tilde r)}}\left[\exp\left(\theta(\tilde r)S(1)\right)I_{\left\lbrace S(1) < \infty\right\rbrace}\right] =p <1.$$
Finally we get \[\mathbb{E}_{(u,\lambda_0)}\left[\exp(-r(X_{S_+(1)}-u))I_{\left\lbrace S_+(1) < \infty\right\rbrace}\right] \leq C \, \sum_{i=1}^\infty p^{i-1} = \frac{C}{1-p} < \infty .\]
\end{proof}
\end{lemma}

\begin{lemma}\label{lem_dri}
The function $\mathbb{P}_{(u,\lambda_0)}\left[\tau_u \leq S_+(1) , \tau_u< \infty\right] e^{Ru}$ is directly Riemann integrable in $u$.
\begin{proof}
Let $r$ be as in Assumption \ref{assumption_r}. Observe that $\alpha(r)<0$ and $\theta(r)>0$ since $r>R>0$. At first, we show that $\mathbb{P}_{(u,\lambda_0)}\left[\tau_u \leq S_+(1), \tau_u < \infty\right]e^{ru}$ is uniformly bounded.
Let $t>0$ be arbitrary but fixed. Then \begin{align*}
\mathbb{P}_{(u,\lambda_0)}\left[ \tau_u \leq (S_+(1) \wedge t)\right] e^{ru}& = \mathbb{E}^{\mathbb{Q}^{(r)}}\left[I_{\left\lbrace \tau_u \leq (S_+(1) \wedge t\right\rbrace} e^{\theta(r) \tau_u} e^{rX_{\tau_u}} e^{\alpha(r) \lambda_{\tau_u}}\right]e^{-\alpha(r) \lambda_0} \\
&\leq 
\mathbb{E}^{\mathbb{Q}^{(r)}}\left[I_{\left\lbrace \tau_u \leq (S_+(1) \wedge t\right\rbrace} e^{\theta(r) \tau_u}\right] e^{-\alpha(r) \lambda_0} \\
&\leq \mathbb{E}^{\mathbb{Q}^{(r)}}\left[ e^{\theta(r) S_+(1)}\right] e^{-\alpha(r) \lambda_0} \\
&= \mathbb{E}_{(u,\lambda_0)} \left[ e^{-rX_{S_+(1)}+ru-\alpha(r) \lambda_{S_+(1)} + \alpha(r) \lambda_0 }I_{\left\lbrace S_+(1) < \infty \right\rbrace}\right]e^{-\alpha(r) \lambda_0} \\
&= \mathbb{E}_{(u,\lambda_0)}\left[I_{\left\lbrace S_+(1) < \infty\right\rbrace}e^{-r(X_{S_+(1)} -u)}\right] e^{-\alpha(r) \lambda_0} < \infty.
\end{align*}
The upper bound is independent of $t$, so by letting $t$ tend to infinity we get \[\mathbb{P}_{(u,\lambda_0)}\left[\tau_u \leq S_+(1), \tau_u < \infty\right]e^{ru} \leq \mathbb{E}_{(u,\lambda_0)}\left[I_{\left\lbrace S_+(1) < \infty\right\rbrace}e^{-r(X_{S_+(1)} -u)}\right] e^{-\alpha(r) \lambda_0}.\] It is even independent of $u$. To see this we consider the Process $R_t = ct- \sum_{i=1}^{N_t} U_i$ and define the random time $T_+(1):= \min\left\lbrace S(i) \, |\, R_{S(i)} < 0 \right\rbrace.$ They are independent of $u$ but under $\mathbb{P}_{(u,\lambda_0)}$ we have almost surely $R_t= X_t-u$ and $T_+(1)= S_+(1)$. By this we see that $X_{S_+(1)}-u = R_{T_+(1)}$ does not depend on $u$.\\
Using the derived boundedness we get that there is some $K>0$ such that \[ \mathbb{P}_{(u,\lambda_0)}\left[\tau_u \leq S_+(1), \tau_u < \infty\right] e^{Ru} \leq K e^{-(r-R)u},\] which is a directly Riemann integrable upper bound. Consequently $$\mathbb{P}_{(u,\lambda_0)}\left[\tau_u \leq S_+(1), \tau_u < \infty\right] e^{Ru}$$ is directly Riemann integrable too.
\end{proof}
\end{lemma}

Let us now focus on the properties of $p(u,x)$.
\begin{lemma}
The function $p(u,x)$ is continuous in $u$ for $u>0$.
\begin{proof}
To prove continuity, we will show that $$\lim_{\varepsilon \to 0} p(u+\varepsilon,x)= \lim_{\varepsilon \to 0} p(u-\varepsilon,x) = p(u,x).$$
We start with the first limit. To do so we will consider a path of our surplus process $X$ with initial capital $u$ and exactly the same path of the process $X^\varepsilon$ with initial capital $u+\varepsilon$. The premium rate $c$, the claim sizes $U_i$ and the counting process $N$ do not depend on the initial capital, hence $X^\varepsilon_t = X_t + \varepsilon$. By the same line of arguments as in the proof of Lemma \ref{lem_dri}, we see that $S_+(1)$ and the condition in the definition of $p$ do not depend on $u$.\\
To be precise, let $\omega \in \Omega$ be an arbitrary event and let us compare the fixed paths of our processes. If $X(\omega)$ gets ruined before $S_+(1)(\omega)$, there is some $\tilde \varepsilon>0$ such that for all $\varepsilon < \tilde \varepsilon$ the path $X^\varepsilon(\omega)$ gets ruined in the same moment. If $X(\omega)$ stays greater or equal to $0$ then $X^\varepsilon$ stays positive for all $\varepsilon>0$. Consequently, we have that $\lim_{\varepsilon \to 0}I_{\left\lbrace \tau_{u+\varepsilon}< S_+(1)\right\rbrace}(\omega) = I_{\left\lbrace \tau_{u}< S_+(1)\right\rbrace}(\omega)$ and by dominated convergence also $p(u+\varepsilon,x) \to p(u,x).$ \\
If we can exclude the case, that $X$ hits exactly the value $0$, then the same arguments hold for $X^{-\varepsilon}_t:= X_t-\varepsilon$.  \\
The infimum of the surplus process can only occur at a jump time of our counting process $N$. Let $T$ be an arbitrary claim time then \begin{align*}
\mathbb{P}_{(u,\lambda_0)} \left[ 
X_T =0 \right] = \mathbb{P}_{(u,\lambda_0)} \left[ X_{T-}-U_{N_T} =0 \right] = 
\mathbb{E}_{(u,\lambda_0)}\left[\mathbb{P}_{(u,\lambda_0)} \left[ X_{T-}-U_{N_T} =0\, \left\vert \, \mathcal{F}_{T-}\right. \right]\right]
\end{align*}
The random variable $U_{N_T}$ is independent of $\mathcal{F}_{T-}$ and its distribution is continuous. Hence the probability of hitting exactly the value $ X_{T-}$ is $0$. Consequently, $\mathbb{P}_{(u,\lambda_0)} \left[  X_T =0 \right]=0$. Since we have only countably many jump times, the event that the surplus process hits $0$ at any jump time has measure $0$ too. 

Hence, $p(u-\varepsilon,x) \to p(u,x).$
Combining these results we get that $p(u,x)$ is continuous in $u$.
\end{proof}
\end{lemma}
\begin{lemma} Under our assumptions $\int_0^u \! p(u,x) e^{Rx} \, B(\mathrm{d}x) $ is directly Riemann integrable.
\begin{proof}
Again let $r$ be as in Assumption \ref{assumption_r}. Then
\begin{align*}
\int_0^u \! p(u,x) e^{Rx} \, B(\mathrm{d}x) 
&\leq e^{Ru} \int_0^u \! p(u,x) \, B(\mathrm{d}x) 
= e^{Ru} \mathbb{P}_{(u,\lambda_0)}\left[\tau_u \leq S_+(1)< \infty\right] \\ 
&\leq  e^{Ru} \mathbb{P}_{(u,\lambda_0)}\left[\tau_u \leq S_+(1), \tau_u< \infty\right] \leq Ke^{-(r-R)u}.
\end{align*}
As before we have a directly Riemann integrable upper bound and therefore $$\int_0^u \! p(u,x) e^{Rx} \, B(\mathrm{d}x)$$ is directly Riemann integrable.
\end{proof} 
\end{lemma}
The continuity of the distribution of $U$ implies that $B$ is not arithmetic. Consequently, we have shown that all conditions of Theorem \ref{renewal_theorem_schmidli} are satisfied. Hence, we can apply it to the renewal equation satisfied by $\psi(u)e^{Ru}$ and obtain our main result.
\begin{theorem}
Under our Assumptions $\lim_{u\to \infty}\psi(u,\lambda_0)e^{Ru}$ exists and is finite. 

\end{theorem}

\section*{Funding}

This research was funded in whole, or in part, by the Austrian Science Fund (FWF) P 33317. For the purpose of open access, the author has applied a CC BY public copyright licence to any Author Accepted Manuscript version arising from this submission.

\end{document}